\documentclass[12pt]{article}
\usepackage{amssymb}
\newcommand{\be}{\begin{equation}}
\newcommand{\ee}{\end{equation}}
\newcommand{\bea}{\begin{eqnarray}}
\newcommand{\eea}{\end{eqnarray}}
\newcommand{\ba}{\begin{array}}
\newcommand{\ea}{\end{array}}

\newcommand{\bc}{\begin{center}}
\newcommand{\ec}{\end{center}}
\newcommand{\ben}{\begin{enumerate}}
\newcommand{\een}{\end{enumerate}}
\newcommand{\bfi}{\begin{figure}}
\newcommand{\efi}{\end{figure}}

\newcommand{\bq}{\begin{quote}}
\newcommand{\eq}{\end{quote}}
\newcommand{\bqu}{\begin{quotation}}
\newcommand{\equ}{\end{quotation}}
\newenvironment{emphit}{\begin{itemize}}{\end{itemize}}
\newcommand{\bemp}{\begin{emphit}}
\newcommand{\eemp}{\end{emphit}}

\newcommand{\bt}{\begin{tabular}}
\newcommand{\et}{\end{tabular}}

\newtheorem{myth}{Theorem}[section]
\newtheorem{mylem}{Lemma}[section]
\newtheorem{mycor}{Corollary}[section]

\newtheorem{mydef}{Definition}[section]

\begin{document}
\date{}
\title{On Smarandache Bryant Schneider Group Of A Smarandache Loop
\footnote{2000 Mathematics Subject Classification. Primary 20NO5 ;
Secondary 08A05.}
\thanks{{\bf Keywords and Phrases : Smarandache Bryant Schneider
group, Smarandache loops, Smarandache $f$,$g$-principal isotopes}}}
\author{T\`em\'it\d{\'o}p\d{\'e} Gb\d{\'o}l\'ah\`an Ja\'iy\'e\d ol\'a\\Department of
Mathematics,\\
Obafemi Awolowo University, Ile Ife,
Nigeria.\\jaiyeolatemitope@yahoo.com, tjayeola@oauife.edu.ng}
\maketitle

\begin{abstract}
The concept of Smarandache Bryant Schneider Group of a Smarandache
loop is introduced. Relationship(s) between the Bryant Schneider
Group and the Smarandache Bryant Schneider Group of an S-loop are
discovered and the later is found to be useful in finding
Smarandache isotopy-isomorphy condition(s) in S-loops just like the
formal is useful in finding isotopy-isomorphy condition(s) in loops.
Some properties of the Bryant Schneider Group of a loop are shown to
be true for the Smarandache Bryant Schneider Group of a Smarandache
loop. Some interesting and useful cardinality formulas are also
established for a type of finite Smarandache loop.
\end{abstract}

\section{Introduction}
\paragraph{}
The study of Smarandache loops was initiated by W. B. Vasantha
Kandasamy in 2002. In her book \cite{phd75}, she defined a
Smarandache loop (S-loop) as a loop with at least a subloop which
forms a subgroup under the binary operation of the loop. For more on
loops and their properties, readers should check \cite{phd3},
\cite{phd41}, \cite{phd39}, \cite{phd49}, \cite{phd42} and
\cite{phd75}. In her book, she introduced over 75 Smarandache
concepts in loops but the concept Smarandache Bryant Schneider Group
which is to be studied here for the first time is not among. In her
first paper \cite{phd83}, she introduced some types of Smarandache
loops. The present author has contributed to the study of
S-quasigroups and S-loops in \cite{sma1}, \cite{sma2} and
\cite{sma3} while Muktibodh \cite{muk} did a study on the first.

Robinson \cite{phd93} introduced the idea of Bryant-Schneider group
of a loop because its importance and motivation stem from the work
of Bryant and Schneider \cite{phd92}. Since the advent of the
Bryant-Schneider group, some studies by Adeniran \cite{phd142},
\cite{ade} and Chiboka \cite{phd94} have been done on it relative to
CC-loops, C-loops and extra loops after Robinson \cite{phd93}
studied the Bryant-Schneider group of a Bol loop. The judicious use
of it was earlier predicted by Robinson \cite{phd93}. As mentioned
in [Section~5, Robinson \cite{phd93}], the Bryant-Schneider group of
a loop is extremely useful in investigating isotopy-isomorphy
condition(s) in loops.

In this study, the concept of Smarandache Bryant Schneider Group of
a Smarandache loop is introduced. Relationship(s) between the Bryant
Schneider Group and the Smarandache Bryant Schneider Group of an
S-loop are discovered and the later is found to be useful in finding
Smarandache isotopy-isomorphy condition(s) in S-loops just like the
formal is useful in finding isotopy-isomorphy condition(s) in loops.
Some properties of the Bryant Schneider Group of a loop are shown to
be true for the Smarandache Bryant Schneider Group of a Smarandache
loop. Some interesting and useful cardinality formulas are also
established for a type of finite Smarandache loop. But first, we
state some important definitions.

\section{Definitions and Notations}
\begin{mydef}
Let $L$ be a non-empty set. Define a binary operation ($\cdot $) on
$L$ : If $x\cdot y\in L~\forall ~x, y\in L$, $(L, \cdot )$ is called
a groupoid. If the system of equations ; $a\cdot x=b$ and $y\cdot
a=b$ have unique solutions for $x$ and $y$ respectively, then $(L,
\cdot )$ is called a quasigroup. Furthermore, if there exists a
unique element $e\in L$ called the identity element such that
$\forall
~x\in L$, $x\cdot e=e\cdot x=x$, $(L, \cdot )$ is called a loop.\\

Furthermore, if there exist at least a non-empty subset $M$ of $L$
such that $(M,\cdot )$ is a non-trivial subgroup of $(L, \cdot )$,
then $L$ is called a Smarandache loop(S-loop) with Smarandache
subgroup(S-subgroup) $M$.
\end{mydef}

The set $SYM(L, \cdot )=SYM(L)$ of all bijections in a loop
$(L,\cdot )$ forms a group called the permutation(symmetric) group
of the loop $(L,\cdot )$. The triple $(U, V, W)$ such that $U, V,
W\in SYM(L, \cdot )$ is called an autotopism of $L$ if and only if
$xU\cdot yV=(x\cdot y)W~\forall ~x, y\in L$. The group of
autotopisms(under componentwise multiplication(\cite{phd3}) of $L$
is denoted by $AUT(L, \cdot )$. If $U=V=W$, then the group
$AUM(L,\cdot )=AUM(L)$ formed by such $U$'s is called the
automorphism group of $(L,\cdot )$. If $L$ is an S-loop with an
arbitrary S-subgroup $H$, then the group $SSYM(L,\cdot )=SSYM(L)$
formed by all $\theta\in SYM(L)$ such that $h\theta\in
H~\forall~h\in H$ is called the Smarandache permutation(symmetric)
group of $L$. Hence, the group $SA(L,\cdot )=SA(L)$ formed by all
$\theta\in SSYM(L)\cap AUM(L)$ is called the Smarandache
automorphism group of $L$.

Let $(G,\cdot )$ be a loop. The bijection $L_x : G\longrightarrow G$
defined as $yL_x=x\cdot y~\forall ~x, y\in G$ is called a left
translation(multiplication) of $G$ while the bijection $R_x :
G\longrightarrow G$ defined as $yR_x=y\cdot x~\forall ~x, y\in G$ is
called a right translation(multiplication) of $G$.

\begin{mydef}\label{1:3}(Robinson \cite{phd93})

Let $(G,\cdot )$ be a loop. A mapping $\theta\in SYM(G,\cdot )$ is a
\textit{special map} for $G$ means that there exist $f,g\in G$ so
that $(\theta R_g^{-1},\theta L_f^{-1},\theta )\in AUT(G,\cdot )$.
\end{mydef}

\begin{mydef}\label{1:4}

Let $(G,\cdot )$ be a Smarandache loop with S-subgroup $(H,\cdot )$.
A mapping $\theta\in SSYM(G,\cdot )$ is a \textit{Smarandache
special map(S-special map)} for $G$ if and only if there exist
$f,g\in H$ such that $(\theta R_g^{-1},\theta L_f^{-1},\theta )\in
AUT(G,\cdot )$.
\end{mydef}

\begin{mydef}\label{1:5}(Robinson \cite{phd93})

Let the set
\begin{displaymath}
BS(G,\cdot )=\{\theta\in SYM(G,\cdot )~:~\exists~f,g\in
G~\ni~(\theta R_g^{-1},\theta L_f^{-1},\theta )\in AUT(G,\cdot )\}
\end{displaymath}
i.e the set of all special maps in a loop, then $BS(G,\cdot )\le
SYM(G,\cdot )$ is called the Bryant-Schneider group of the loop
$(G,\cdot )$.
\end{mydef}

\begin{mydef}\label{1:6}

Let the set
\begin{displaymath}
SBS(G,\cdot )=\{\theta\in SSYM(G,\cdot )~:~\textrm{there
exist}~f,g\in H~\ni~(\theta R_g^{-1},\theta L_f^{-1},\theta )\in
AUT(G,\cdot )\}
\end{displaymath}
i.e the set of all S-special maps in a S-loop, then $SBS(G,\cdot )$
is called the Smarandache Bryant-Schneider group(SBS group) of the
S-loop $(G,\cdot )$ with S-subgroup $H$ if $SBS(G,\cdot )\le
SYM(G,\cdot )$.
\end{mydef}

\begin{mydef}\label{1:7}
The triple $\phi =(R_g,L_f,I)$ is called an $f,g$-principal
isotopism of a loop $(G,\cdot )$ onto a loop $(G,\circ )$ if and
only if
\begin{displaymath}
x\cdot y=xR_g\circ yL_f~\forall~x,y\in G~\textrm{or}~x\circ
y=xR_g^{-1}\cdot yL_f^{-1}~\forall~x,y\in G.
\end{displaymath}
$f$ and $g$ are called translation elements of $G$ or at times
written in the pair form $(g,f)$, while $(G,\circ )$ is called an
$f,g$-principal isotope of $(G,\cdot )$.\\

On the other hand, $(G,\otimes )$ is called a Smarandache
$f,g$-principal isotope of $(G,\oplus )$ if for some $f,g\in S$,
\begin{displaymath}
xR_g\otimes yL_f=(x\oplus y)~\forall~x,y\in G
\end{displaymath}
where $(S,\oplus )$ is a S-subgroup of $(G,\oplus )$. In these
cases, $f$ and $g$ are called Smarandache elements(S-elements).

Let $(L, \cdot )$ and $(G, \circ )$ be S-loops with S-subgroups $L'$
and $G'$ respectively such that $xA\in G'~\forall~x\in L'$, where
$A~:~(L, \cdot )\longrightarrow (G, \circ )$. Then the mapping $A$
is called a Smarandache isomorphism if $(L, \cdot )\cong (G, \circ
)$, hence we write $(L, \cdot )\succsim (G, \circ )$. An S-loop $(L,
\cdot )$ is called a G-Smarandache loop(GS-loop) if and only if $(L,
\cdot )\succsim (G, \circ )$ for all S-loop isotopes $(G, \circ )$
of $(L, \cdot )$.
\end{mydef}

\begin{mydef}\label{1:7.1}

Let $(G,\cdot )$ be a Smarandache loop with an S-subgroup $H$.
\begin{displaymath}
\Omega (G,\cdot )=\bigg\{(\theta R_g^{-1},\theta L_f^{-1},\theta
)\in AUT(G,\cdot )~\textrm{for some}~f,g\in H~:~h\theta \in
H~\forall~h\in H\bigg\}
\end{displaymath}
\end{mydef}

\section{Main Results}
\subsection*{The Smarandache Bryant Schneider Group}
\begin{myth}\label{1:10}
Let $(G,\cdot )$ be a Smarandache loop. $SBS(G,\cdot )\le BS(G,\cdot
)$.
\end{myth}
{\bf Proof}\\
Let $(G,\cdot )$ be an S-loop with S-subgroup $H$. Comparing
Definition~\ref{1:5} and Definition~\ref{1:6}, it can easily be
observed that $SBS(G,\cdot )\subset BS(G,\cdot )$. The case
$SBS(G,\cdot )\subseteq BS(G,\cdot )$ is possible when $G=H$ where
$H$ is the S-subgroup of $G$ but this will be a contradiction since
$G$ is an S-loop.
\begin{description}
\item[Identity]
If $I$ is the identity mapping on $G$, then $hI=h\in H~\forall~h\in
H$ and there exists $e\in H$ where $e$ is the identity element in
$G$ such that $(IR_e^{-1},IL_e^{-1},I)=(I,I,I)\in AUT(G,\cdot)$. So,
$I\in SBS(G,\cdot )$. Thus $SBS(G,\cdot )$ is non-empty.
\item[Closure and Inverse]
Let $\alpha,\beta\in SBS(G,\cdot )$. Then there exist
$f_1,g_1,f_2,g_2\in H$ such that
\begin{displaymath}
A=(\alpha R_{g_1}^{-1},\alpha L_{f_1}^{-1},\alpha ),~B=(\beta
R_{g_2}^{-1},\beta L_{f_2}^{-1},\beta )\in AUT(G,\cdot ).
\end{displaymath}
\begin{displaymath}
AB^{-1}=(\alpha R_{g_1}^{-1},\alpha L_{f_1}^{-1},\alpha
)(R_{g_2}\beta^{-1},L_{f_2}\beta^{-1},\beta^{-1})
\end{displaymath}
\begin{displaymath}
=(\alpha R_{g_1}^{-1}R_{g_2}\beta^{-1},\alpha
L_{f_1}^{-1}L_{f_2}\beta^{-1},\alpha\beta^{-1} )\in AUT(G,\cdot ).
\end{displaymath}
Let $\delta =\beta R_{g_1}^{-1}R_{g_2}\beta^{-1}$ and $\gamma=\beta
L_{f_1}^{-1}L_{f_2}\beta^{-1}$. Then,
\begin{displaymath}
(\alpha \beta^{-1}\delta,\alpha \beta^{-1}\gamma,\alpha\beta^{-1}
)\in AUT(G,\cdot )\Leftrightarrow (x\alpha \beta^{-1}\delta)\cdot
(y\alpha \beta^{-1}\gamma)=(x\cdot y)\alpha\beta^{-1}~\forall~x,y\in
G.
\end{displaymath}
Putting $y=e$ and replacing $x$ by $x\beta\alpha^{-1}$, we have
$(x\delta)\cdot (e\alpha \beta^{-1}\gamma)=x$ for all $x\in G$.
Similarly, putting $x=e$ and replacing $y$ by $y\beta\alpha^{-1}$,
we have $(e\alpha \beta^{-1}\delta)\cdot (y\gamma)=y$ for all $y\in
G$. Thence, $x\delta R_{(e\alpha \beta^{-1}\gamma)}=x$ and $y\gamma
L_{(e\alpha \beta^{-1}\delta)}=y$ which implies that
\begin{displaymath}
\delta= R_{(e\alpha \beta^{-1}\gamma )}^{-1}~\textrm{and}~\gamma
=L_{(e\alpha \beta^{-1}\delta)}^{-1}.
\end{displaymath}
Thus, since $g=e\alpha \beta^{-1}\gamma,~f=e\alpha
\beta^{-1}\delta\in H$ then
\begin{displaymath}
AB^{-1}=(\alpha \beta^{-1}R_g^{-1},\alpha
\beta^{-1}L_f^{-1},\alpha\beta^{-1} )\in AUT(G,\cdot
)\Leftrightarrow \alpha\beta^{-1}\in SBS(G,\cdot ).
\end{displaymath}
\end{description}
$\therefore SBS(G,\cdot )\le BS(G,\cdot )$.

\begin{mycor}\label{1:11} Let $(G,\cdot
)$ be a Smarandache loop. Then, $SBS(G,\cdot )\le SSYM(G,\cdot )\le
SYM(G,\cdot )$. Hence, $SBS(G,\cdot )$ is the Smarandache
Bryant-Schneider group(SBS group) of the S-loop $(G,\cdot )$.
\end{mycor}
{\bf Proof}\\
Although the fact that $SBS(G,\cdot )\le SYM(G,\cdot )$ follows from
Theorem~\ref{1:10} and the fact in [Theorem~1, \cite{phd93}] that
$BS(G,\cdot )\le SYM(G,\cdot )$. Nevertheless, it can also be traced
from the facts that $SBS(G,\cdot )\le SSYM(G,\cdot )$ and
$SSYM(G,\cdot )\le SYM(G,\cdot )$.

It is easy to see that $SSYM(G,\cdot )\subset SYM(G,\cdot )$ and
that $SBS(G,\cdot )\subset SSYM(G,\cdot )$ while the trivial cases
$SSYM(G,\cdot )\subseteq SYM(G,\cdot )$ and $SBS(G,\cdot )\subseteq
SSYM(G,\cdot )$ will contradict the fact that $G$ is an S-loop
because these two are possible if the S-subgroup $H$ is $G$.
Reasoning through the axioms of a group, it is easy to show that
$SSYM(G,\cdot )\le SYM(G,\cdot )$. By using the same steps in
Theorem~\ref{1:10}, it will be seen that $SBS(G,\cdot )\le
SSYM(G,\cdot )$.

\subsection*{The SBS Group of a Smarandache $f,g$-principal isotope}
\begin{myth}\label{1:12}
Let $(G,\cdot )$ be a S-loop with a Smarandache $f,g$-principal
isotope $(G,\circ )$. Then, $(G,\circ )$ is an S-loop.
\end{myth}
{\bf Proof}\\
Let $(G,\cdot )$ be an S-loop, then there exist an S-subgroup
$(H,\cdot )$ of $G$. If $(G,\circ )$ is a Smarandache
$f,g$-principal isotope of $(G,\cdot )$, then
\begin{displaymath}
x\cdot y=xR_g\circ yL_f~\forall~x,y\in G~\textrm{which
implies}~x\circ y=xR_g^{-1}\cdot yL_f^{-1}~\forall~x,y\in G
\end{displaymath}
where $f,g\in H$. So
\begin{displaymath}
h_1\circ h_2=h_1R_g^{-1}\cdot h_2L_f^{-1}~\forall~h_1,h_2\in
H~\textrm{for some}~f,g\in H.
\end{displaymath}
Let us now consider the set $H$ under the operation "$\circ$". That
is the pair $(H,\circ )$.
\begin{description}
\item[Groupoid] Since $f,g\in H$, then by the definition $h_1\circ h_2=h_1R_g^{-1}\cdot
h_2L_f^{-1},~h_1\circ h_2\in H~\forall~h_1,h_2\in H$ since $(H,\cdot
)$ is a groupoid. Thus, $(H,\circ )$ is a groupoid.
\item[Quasigroup] With the definition $h_1\circ h_2=h_1R_g^{-1}\cdot
h_2L_f^{-1}~\forall~h_1,h_2\in H$, it is clear that $(H,\circ )$ is
a quasigroup since $(H,\cdot )$ is a quasigroup.
\item[Loop] It can easily be seen that $f\cdot g$ is an identity element
in $(H,\circ )$. So, $(H,\circ )$ is a loop.
\item[Group] Since $(H,\cdot )$ is a associative, it is easy to show that $(H,\circ )$ is
associative.
\end{description}
Hence, $(H,\circ )$ is an S-subgroup in $(G,\circ )$ since the
latter is a loop(a quasigroup with identity element $f\cdot g$).
Therefore, $(G,\circ )$ is an S-loop.

\begin{myth}\label{1:8}
Let $(G,\cdot )$ be a Smarandache loop with an S-subgroup $(H,\cdot
)$. A mapping $\theta\in SYM(G,\cdot )$ is a S-special map if and
only if $\theta$ is an S-isomorphism of $(G,\cdot )$ onto some
Smarandache $f,g$-principal isotopes $(G,\circ )$ where $f,g\in H$.
\end{myth}
{\bf Proof}\\
By Definition~\ref{1:4}, a mapping $\theta\in SSYM(G)$ is a
S-special map implies there exist $f,g\in H$ such that $(\theta
R_g^{-1},\theta L_f^{-1},\theta )\in AUT(G,\cdot )$. It can be
observed that
\begin{displaymath}
(\theta R_g^{-1},\theta L_f^{-1},\theta )=(\theta ,\theta ,\theta
)(R_g^{-1},L_f^{-1},I)\in AUT(G,\cdot ).
\end{displaymath}
But since $(R_g^{-1},L_f^{-1},I)~:~(G,\circ )\longrightarrow
(G,\cdot )$ then for $(\theta R_g^{-1},\theta L_f^{-1},\theta )\in
AUT(G,\cdot)$ we must have $(\theta ,\theta ,\theta )~:~(G,\cdot
)\longrightarrow (G,\circ )$ which means $(G,\cdot )\cong
^{^\theta}(G,\circ )$, hence $(G,\cdot )\succsim^{^\theta}(G,\circ
)$ because $(H,\cdot )\theta =(H,\circ )$. $(R_g,L_f,I)~:~(G,\cdot
)\longrightarrow (G,\circ )$ is an $f,g$-principal isotopism so
$(G,\circ )$ is a Smarandache $f,g$-principal isotope of $(G,\cdot
)$ by Theorem~\ref{1:12}.

Conversely, if $\theta$ is an S-isomorphism of $(G,\cdot )$ onto
some Smarandache $f,g$-principal isotopes $(G,\circ )$ where $f,g\in
H$ such that $(H,\cdot )$ is a S-subgroup of $(G,\cdot )$ means
$(\theta ,\theta ,\theta )~:~(G,\cdot )\longrightarrow (G,\circ )$,
$(R_g,L_f,I)~:~(G,\cdot )\longrightarrow (G,\circ )$ which implies
$(R_g^{-1},L_f^{-1},I)~:~(G,\circ )\longrightarrow (G,\cdot )$ and
$(H,\cdot )\theta =(H,\circ )$. Thus, $(\theta R_g^{-1},\theta
L_f^{-1},\theta )\in AUT(G,\cdot)$. Therefore, $\theta$ is a
S-special map because $f,g\in H$.

\begin{mycor}\label{1:9}
Let $(G,\cdot )$ be a Smarandache loop with a an S-subgroup
$(H,\cdot )$. A mapping $\theta\in SBS(G,\cdot )$ if and only if
$\theta$ is an S-isomorphism of $(G,\cdot )$ onto some Smarandache
$f,g$-principal isotopes $(G,\circ )$ such that $f,g\in H$ where
$(H,\cdot )$ is an S-subgroup of $(G,\cdot )$.
\end{mycor}
{\bf Proof}\\
This follows from Definition~\ref{1:6} and Theorem~\ref{1:8}.

\begin{myth}\label{1:12.1}
Let $(G,\cdot )$ and $(G,\circ )$ be S-loops. $(G,\circ )$ is a
Smarandache $f,g$-principal isotope of $(G,\cdot )$ if and only if
$(G,\cdot )$ is a Smarandache $g,f$-principal isotope of $(G,\circ
)$.
\end{myth}
{\bf Proof}\\
Let $(G,\cdot )$ and $(G,\circ )$ be S-loops such that if $(H,\cdot
)$ is an S-subgroup in $(G,\cdot )$, then $(H,\circ )$ is an
S-subgroup of $(G,\circ )$. The left and right translation maps
relative to an element $x$ in $(G,\circ )$ shall be denoted by
${\cal L}_x$ and ${\cal R}_x$ respectively.

If $(G,\circ )$ is a Smarandache $f,g$-principal isotope of
$(G,\cdot )$ then, $x\cdot y=xR_g\circ yL_f~\forall~x,y\in G$ for
some $f,g\in H$. Thus, $xR_y=xR_g{\cal R }_{yL_f}$ and
$yL_x=yL_f{\cal L}_{xR_g}~x,y\in G$ and we have $R_y=R_g{\cal R
}_{yL_f}$ and $L_x=L_f{\cal L}_{xR_g}~x,y\in G$. So, ${\cal R
}_y=R_g^{-1}R_{yL_f^{-1}}$ and ${\cal
L}_{x}=L_f^{-1}L_{xR_g^{-1}}=~x,y\in G$. Putting $y=f$ and $x=g$
respectively, we now get ${\cal R
}_f=R_g^{-1}R_{fL_f^{-1}}=R_g^{-1}$ and ${\cal
L}_g=L_f^{-1}L_{gR_g^{-1}}=L_f^{-1}$. That is, ${\cal R
}_f=R_g^{-1}$ and ${\cal L}_g=L_f^{-1}$ for some $f,g\in H$.

Recall that
\begin{displaymath}
x\cdot y=xR_g\circ yL_f~\forall~x,y\in G~\Leftrightarrow~x\circ
y=xR_g^{-1}\cdot yL_f^{-1}~\forall~x,y\in G.
\end{displaymath}
So using the last two translation equations,
\begin{displaymath}
x\circ y=x{\cal R}_f\cdot y{\cal L}_g~\forall~x,y\in
G~\Leftrightarrow~\textrm{the triple}~({\cal R}_f,{\cal L}_g,I
)~:~(G,\circ )\longrightarrow (G,\cdot )
\end{displaymath}
is a Smarandache $g,f$-principal isotopism. Therefore, $(G,\cdot )$
is a Smarandache $g,f$-principal isotope of $(G,\circ )$.\\

The proof of the converse is achieved by doing the reverse of the
procedure described above.

\begin{myth}\label{1:13}
If $(G,\cdot )$ is an S-loop with a Smarandache $f,g$-principal
isotope $(G,\circ )$, then $SBS(G,\cdot )=SBS(G,\circ )$.
\end{myth}
{\bf Proof}\\
Let $(G,\circ )$ be the Smarandache $f,g$-principal isotope of the
S-loop $(G,\cdot )$ with S-subgroup $(H,\cdot )$. By
Theorem~\ref{1:12}, $(G,\circ )$ is an S-loop with S-subgroup
$(H,\circ )$. The left and right translation maps relative to an
element $x$ in $(G,\circ )$ shall be denoted by ${\cal L}_x$ and
${\cal R}_x$ respectively.

Let $\alpha\in SBS(G,\cdot )$, then there exist $f_1,g_1\in H$ so
that $(\alpha R_{g_1}^{-1},\alpha L_{f_1}^{-1},\alpha )\in
AUT(G,\cdot)$. Recall that the triple $(R_{g_1},L_{f_1},I
)~:~(G,\cdot )\longrightarrow (G,\circ )$ is a Smarandache
$f,g$-principal isotopism, so $x\cdot y=xR_g\circ
yL_f~\forall~x,y\in G$ and this implies
\begin{displaymath}
R_x=R_g{\cal R}_{xL_f}~\textrm{and}~L_x=L_f{\cal
L}_{xR_g}~\forall~x\in G~\textrm{which also implies that}
\end{displaymath}
\begin{displaymath}
{\cal R}_{xL_f}=R_g^{-1}R_x~\textrm{and}~{\cal
L}_{xR_g}=L_f^{-1}L_x~\forall~x\in G~\textrm{which finally gives}
\end{displaymath}
\begin{displaymath}
{\cal R}_x=R_g^{-1}R_{xL_f^{-1}}~\textrm{and}~{\cal
L}_x=L_f^{-1}L_{xR_g^{-1}}~\forall~x\in G.
\end{displaymath}
Set $f_2=f\alpha R_{g_1}^{-1}R_g$ and $g_2=g\alpha L_{f_1}^{-1}L_f$.
Then
\begin{equation}\label{eq:1}
{\cal R}_{g_2}=R_g^{-1}R_{g\alpha
L_{f_1}^{-1}L_fL_f^{-1}}=R_g^{-1}R_{g\alpha L_{f_1}^{-1}}
\end{equation}
\begin{equation}\label{eq:1.2}
~\textrm{and}~{\cal L}_{f_2}=L_f^{-1}L_{f\alpha
R_{g_1}^{-1}R_gR_g^{-1}}=L_f^{-1}L_{f\alpha
R_{g_1}^{-1}}~\forall~x\in G.
\end{equation}
Since, $(\alpha R_{g_1}^{-1},\alpha L_{f_1}^{-1},\alpha )\in
AUT(G,\cdot)$, then
\begin{equation}\label{eq:2}
(x\alpha R_{g _1}^{-1})\cdot (y\alpha L_{f_1}^{-1})=(x\cdot
y)\alpha~\forall~x,y\in G.
\end{equation}
Putting $y=g$ and $x=f$ separately in the last equation,
\begin{displaymath}
x\alpha R_{g_1}^{-1}R_{(g\alpha
L_{f_1}^{-1})}=xR_g\alpha~\textrm{and}~y\alpha
L_{f_1}^{-1}L_{(f\alpha R_{g _1}^{-1})}=yL_f\alpha~\forall~x,y\in G.
\end{displaymath}
Thus by applying (\ref{eq:1}) and (\ref{eq:1.2}), we now have
\begin{equation}\label{eq:3}
\alpha R_{g _1}^{-1}=R_g\alpha R_{(g\alpha
L_{f_1}^{-1})}^{-1}=R_g\alpha {\cal R}_{g_2}^{-1}R_g^{-1}
~\textrm{and}~\alpha L_{f_1}^{-1}=L_f\alpha L_{(f\alpha R_{g
_1}^{-1})}^{-1}=L_f\alpha {\cal L}_{f_2}^{-1}L_f^{-1}.
\end{equation}
We shall now compute $(x\circ y)\alpha$ by (\ref{eq:2}) and
(\ref{eq:3}) and then see the outcome.

$(x\circ y)\alpha =(xR_g^{-1}\cdot yL_f^{-1})\alpha =xR_g^{-1}\alpha
R_{g _1}^{-1}\cdot yL_f^{-1}\alpha L_{f_1}^{-1}=xR_g^{-1}R_g\alpha
{\cal R}_{g_2}^{-1}R_g^{-1}\cdot yL_f^{-1}L_f\alpha {\cal
L}_{f_2}^{-1}L_f^{-1}=x\alpha {\cal R}_{g_2}^{-1}R_g^{-1}\cdot
y\alpha {\cal L}_{f_2}^{-1}L_f^{-1}=x\alpha {\cal R}_{g_2}^{-1}\circ
y\alpha {\cal L}_{f_2}^{-1}~\forall~x,y\in G$.

Thus,
\begin{displaymath}
(x\circ y)\alpha =x\alpha {\cal R}_{g_2}^{-1}\circ y\alpha {\cal
L}_{f_2}^{-1}~\forall~x,y\in G\Leftrightarrow (\alpha {\cal
R}_{g_2}^{-1},\alpha {\cal L}_{f_2}^{-1},\alpha )\in
AUT(G,\circ)\Leftrightarrow \alpha\in SBS(G,\circ ).
\end{displaymath}
Whence, $SBS(G,\cdot )\subseteq SBS(G,\circ )$.\\

Since $(G,\circ )$ is the Smarandache $f,g$-principal isotope of the
S-loop $(G,\cdot )$, then by Theorem~\ref{1:12.1}, $(G,\cdot )$ is
the Smarandache $g,f$-principal isotope of $(G,\circ )$. So
following the steps above, it can similarly be shown that
$SBS(G,\circ )\subseteq SBS(G,\cdot )$. Therefore, the conclusion
that $SBS(G,\cdot )=SBS(G,\circ )$ follows.

\subsection*{Cardinality Formulas}
\begin{myth}\label{1:14}
Let $(G,\cdot )$ be a finite Smarandache loop with $n$ distinct
S-subgroups. If the SBS group of $(G,\cdot )$ relative to an
S-subgroup $(H_i,\cdot )$ is denoted by $SBS_i(G,\cdot )$, then
\begin{displaymath}
|BS(G,\cdot )|=\frac{1}{n}\sum_{i=1}^n|SBS_i(G,\cdot )|~[BS(G,\cdot
):SBS_i(G,\cdot )].
\end{displaymath}
\end{myth}
{\bf Proof}\\
Let the $n$ distinct S-subgroups of $G$ be denoted by
$H_i,~i=1,2,\cdots n$. Note here that $H_i\ne
H_j~\forall~i,j=1,2,\cdots n$. By Theorem~\ref{1:10}, $SBS_i(G,\cdot
)\le BS(G,\cdot )~\forall~i=1,2,\cdots n$. Hence, by the Lagrange's
theorem of classical group theory,
\begin{displaymath}
|BS(G,\cdot )|=|SBS_i(G,\cdot )|~[BS(G,\cdot ):SBS_i(G,\cdot
)]~\forall~i=1,2,\cdots n.
\end{displaymath}
Thus, adding the equation above for all $i=1,2,\cdots n$, we get
\begin{displaymath}
n|BS(G,\cdot )|=\sum_{i=1}^n|SBS_i(G,\cdot )|~[BS(G,\cdot
):SBS_i(G,\cdot )]~\forall~i=1,2,\cdots n,~\textrm{thence},
\end{displaymath}
\begin{displaymath}
|BS(G,\cdot )|=\frac{1}{n}\sum_{i=1}^n|SBS_i(G,\cdot )|~[BS(G,\cdot
):SBS_i(G,\cdot )].
\end{displaymath}

\begin{myth}\label{1:15}
Let $(G,\cdot )$ be a Smarandache loop. Then, $\Omega (G,\cdot )\le
AUT(G,\cdot)$.
\end{myth}
{\bf Proof}\\
Let $(G,\cdot )$ be an S-loop with S-subgroup $H$. By
Definition~\ref{1:7.1}, it can easily be observed that $\Omega
(G,\cdot )\subseteq AUT(G,\cdot)$.
\begin{description}
\item[Identity]
If $I$ is the identity mapping on $G$, then $hI=h\in H~\forall~h\in
H$ and there exists $e\in H$ where $e$ is the identity element in
$G$ such that $(IR_e^{-1},IL_e^{-1},I)=(I,I,I)\in AUT(G,\cdot)$. So,
$(I,I,I)\in \Omega (G,\cdot )$. Thus $\Omega (G,\cdot )$ is
non-empty.
\item[Closure and Inverse]
Let $A,B\in \Omega(G,\cdot )$. Then there exist $\alpha, \beta\in
SSYM (G,\cdot)$ and some $f_1,g_1,f_2,g_2\in H$ such that
\begin{displaymath}
A=(\alpha R_{g_1}^{-1},\alpha L_{f_1}^{-1},\alpha ),~B=(\beta
R_{g_2}^{-1},\beta L_{f_2}^{-1},\beta )\in AUT(G,\cdot ).
\end{displaymath}
\begin{displaymath}
AB^{-1}=(\alpha R_{g_1}^{-1},\alpha L_{f_1}^{-1},\alpha
)(R_{g_2}\beta^{-1},L_{f_2}\beta^{-1},\beta^{-1})
\end{displaymath}
\begin{displaymath}
=(\alpha R_{g_1}^{-1}R_{g_2}\beta^{-1},\alpha
L_{f_1}^{-1}L_{f_2}\beta^{-1},\alpha\beta^{-1} )\in AUT(G,\cdot ).
\end{displaymath}
Using the same techniques for the proof of closure and inverse in
Theorem~\ref{1:10} here and by letting $\delta =\beta
R_{g_1}^{-1}R_{g_2}\beta^{-1}$ and $\gamma=\beta
L_{f_1}^{-1}L_{f_2}\beta^{-1}$, it can be shown that,
\begin{displaymath}
AB^{-1}=(\alpha \beta^{-1}R_g^{-1},\alpha
\beta^{-1}L_f^{-1},\alpha\beta^{-1} )\in AUT(G,\cdot
)~\textrm{where}~g=e\alpha \beta^{-1}\gamma,~f=e\alpha
\beta^{-1}\delta\in H
\end{displaymath}
\begin{displaymath}
\textrm{such that}~\alpha \beta^{-1}\in SSYM(G,\cdot
)\Leftrightarrow AB^{-1}\in \Omega (G,\cdot ).
\end{displaymath}
\end{description}
$\therefore \Omega (G,\cdot )\le AUT(G,\cdot)$.

\begin{myth}\label{1:16}
Let $(G,\cdot )$ be a Smarandache loop with an S-subgroup $H$ such
that $f,g\in H$ and $\alpha \in SBS(G,\cdot)$. If the mapping
\begin{displaymath}
\Phi~:~\Omega (G,\cdot )\longrightarrow SBS(G,\cdot)~~\textrm{is
defined as} ~~\Phi~:~(\alpha R_g^{-1},\alpha L_f^{-1},\alpha
)\mapsto \alpha,
\end{displaymath}
then $\Phi$ is an homomorphism.
\end{myth}
{\bf Proof}\\
Let $A,B\in \Omega(G,\cdot )$. Then there exist $\alpha, \beta\in
SSYM (G,\cdot)$ and some $f_1,g_1,f_2,g_2\in H$ such that
\begin{displaymath}
A=(\alpha R_{g_1}^{-1},\alpha L_{f_1}^{-1},\alpha ),~B=(\beta
R_{g_2}^{-1},\beta L_{f_2}^{-1},\beta )\in AUT(G,\cdot ).
\end{displaymath}
$\Phi (AB)=\Phi [(\alpha R_{g_1}^{-1},\alpha L_{f_1}^{-1},\alpha
)(\beta R_{g_2}^{-1},\beta L_{f_2}^{-1},\beta )]=\Phi (\alpha
R_{g_1}^{-1}\beta R_{g_2}^{-1},\alpha L_{f_1}^{-1}\beta
L_{f_2}^{-1},\alpha\beta )$. It will be good if this can be written
as; $\Phi (AB)=\Phi (\alpha\beta\delta ,\alpha
\beta\gamma,\alpha\beta)$ such that $h\alpha\beta\in H~\forall~h\in
H$ and $\delta =R_g^{-1}$, $\gamma =L_f^{-1}$ for some $g,f\in H$.

This is done as follows: If
\begin{displaymath}
(\alpha R_{g_1}^{-1}\beta R_{g_2}^{-1},\alpha L_{f_1}^{-1}\beta
L_{f_2}^{-1},\alpha\beta )=(\alpha\beta\delta ,\alpha
\beta\gamma,\alpha\beta)\in AUT(G,\cdot )~\textrm{then},
\end{displaymath}
\begin{displaymath}
x\alpha\beta\delta\cdot y\alpha \beta\gamma =(x\cdot
y)\alpha\beta~\forall~x,y\in G.
\end{displaymath}
Put $y=e$ and replace $x$ by $x\beta^{-1}\alpha^{-1}$ then
$x\delta\cdot e\alpha \beta\gamma =x\Leftrightarrow \delta
=R_{e\alpha \beta\gamma}^{-1}$.

Similarly, put $x=e$ and replace $y$ by $y\beta^{-1}\alpha^{-1}$.
Then, $e\alpha\beta\delta\cdot y\gamma =y\Leftrightarrow \gamma =
L_{e\alpha\beta\delta}^{-1}$. So,
\begin{displaymath}
\Phi (AB)=(\alpha\beta R_{e\alpha \beta\gamma}^{-1} ,\alpha \beta
L_{e\alpha\beta\delta}^{-1},\alpha\beta)=\alpha\beta =\Phi(\alpha
R_{g_1}^{-1},\alpha L_{f_1}^{-1},\alpha )\Phi (\beta
R_{g_2}^{-1},\beta L_{f_2}^{-1},\beta )=\Phi (A)\Phi (B).
\end{displaymath}
$\therefore \Phi$ is an homomorphism.

\begin{myth}\label{1:17}
Let $(G,\cdot )$ be a Smarandache loop with an S-subgroup $H$ such
that $f,g\in H$ and $\alpha \in SSYM(G,\cdot)$. If the mapping
\begin{displaymath}
\Phi~:~\Omega (G,\cdot )\longrightarrow SBS(G,\cdot)~~\textrm{is
defined as} ~~\Phi~:~(\alpha R_g^{-1},\alpha L_f^{-1},\alpha
)\mapsto \alpha
\end{displaymath}
then,
\begin{displaymath} A=(\alpha R_g^{-1},\alpha
L_f^{-1},\alpha )\in \ker\Phi~\textrm{if and only if}~\alpha
\end{displaymath}
is the identity map on $G$, $g\cdot f$ is the identity element of
$(G,\cdot )$ and $g\in N_\mu (G,\cdot )$ the middle nucleus of
$(G,\cdot )$.
\end{myth}
{\bf Proof}\\
\begin{description}
\item[Necessity] $\ker \Phi =\{A\in \Omega(G,\cdot ):\Phi (A)=I\}$. So, if $A=(\alpha R_{g_1}^{-1},\alpha L_{f_1}^{-1},\alpha )\in
\ker\Phi$, then $\Phi (A)=\alpha =I$. Thus, $A=(R_{g_1}^{-1},
L_{f_1}^{-1},I  )\in AUT(G,\cdot)\Leftrightarrow$
\begin{equation}\label{eq:4}
x\cdot y=xR_g^{-1}\cdot yL_f^{-1}~\forall~x,y\in G.
\end{equation}
Replace $x$ by $xR_g$ and $y$ by $yL_f$ in (\ref{eq:4}) to get
\begin{equation}\label{eq:5}
x\cdot y=xg\cdot fy~\forall~x,y\in G.
\end{equation}
Putting $x=y=e$ in (\ref{eq:5}), we get $g\cdot f=e$. Replace $y$ by
$yL_f^{-1}$ in (\ref{eq:5}) to get
\begin{equation}\label{eq:6}
x\cdot yL_f^{-1}=xg\cdot y~\forall~x,y\in G.
\end{equation}
Put $x=e$ in (\ref{eq:6}), then we have $yL_f^{-1}=g\cdot
y~\forall~y\in G$ and  so (\ref{eq:6}) now becomes
\begin{displaymath}
x\cdot (gy)=xg\cdot y~\forall~x,y\in G\Leftrightarrow g\in N_\mu
(G,\cdot ).
\end{displaymath}
\item[Sufficiency] Let $\alpha$ be the identity map on $G$, $g\cdot f$ the identity element of
$(G,\cdot )$ and $g\in N_\mu (G,\cdot )$. Thus, $fg\cdot f=f\cdot
gf=fe=f$. Thus, $f\cdot g=e$. Then also, $y=fg\cdot y=f\cdot
gy~\forall~y\in G$ which results into $yL_f^{-1}=gy~\forall~y\in G$.
Thus, it can be seen that $x\alpha R_g^{-1}\cdot y\alpha
L_f^{-1}=xR_g^{-1}\cdot yL_f^{-1}=xR_g^{-1}\alpha\cdot
yL_f^{-1}\alpha=xR_g^{-1}\cdot yL_f^{-1}=xR_g^{-1}\cdot
gy=(xR_g^{-1}\cdot g)y=xR_g^{-1}R_g\cdot y=x\cdot y=(x\cdot
y)\alpha~\forall~x,y\in G$. Thus, $\Phi (A)=\Phi (\alpha
R_g^{-1},\alpha L_f^{-1},\alpha )=\Phi (R_g^{-1},L_f^{-1},I
)=I\Rightarrow A\in \ker\Phi$.
\end{description}

\begin{myth}\label{1:18}
Let $(G,\cdot )$ be a Smarandache loop with an S-subgroup $H$ such
that $f,g\in H$ and $\alpha \in SSYM(G,\cdot)$. If the mapping
\begin{displaymath}
\Phi~:~\Omega (G,\cdot )\longrightarrow SBS(G,\cdot)~~\textrm{is
defined as} ~~\Phi~:~(\alpha R_g^{-1},\alpha L_f^{-1},\alpha
)\mapsto \alpha
\end{displaymath}
then,
\begin{displaymath}
|N_\mu (G,\cdot )|=|\ker \Phi |~\textrm{and}~|\Omega (G,\cdot
)|=|SBS(G,\cdot)||N_\mu (G,\cdot )|.
\end{displaymath}
\end{myth}
{\bf Proof}\\
Let the identity map on $G$ be $I$. Using Theorem~\ref{1:17}, if
\begin{displaymath}
g\theta =(R_g^{-1},L_{g^{-1}}^{-1},I )~\forall~g\in N_\mu (G,\cdot
)~\textrm{then},~\theta~:~N_\mu (G,\cdot )\longrightarrow \ker \Phi.
\end{displaymath}
$\theta$ is easily seen to be a bijection, hence $|N_\mu (G,\cdot
)|=|\ker \Phi |$.

Since $\Phi$ is an homomorphism by Theorem~\ref{1:16}, then by the
first isomorphism theorem in classical group theory, $\Omega
(G,\cdot )/\ker \Phi\cong \textrm{Im}\Phi$. $\Phi$ is clearly onto,
so $\textrm{Im}\Phi = SBS(G,\cdot)$, so that $\Omega (G,\cdot )/\ker
\Phi\cong SBS(G,\cdot)$. Thus, $|\Omega (G,\cdot )/\ker
\Phi|=|SBS(G,\cdot)|$. By Lagrange's theorem, $|\Omega (G,\cdot
)|=|\ker \Phi||\Omega (G,\cdot )/\ker \Phi|$, so, $|\Omega (G,\cdot
)|=|\ker \Phi||SBS(G,\cdot)|$, $\therefore |\Omega (G,\cdot
)|=|N_\mu (G,\cdot )||SBS(G,\cdot)|$.

\begin{myth}\label{1:19}
Let $(G,\cdot )$ be a Smarandache loop with an S-subgroup $H$. If
\begin{displaymath}
\Theta (G,\cdot )=\Big\{(f,g)\in H\times H~:~(G,\circ )\succsim
(G,\cdot )
\end{displaymath}
\begin{displaymath}
\textrm{for}~(G,\circ )~\textrm{the Smarandache principal
$f,g$-isotope of}~(G,\cdot )\Big\}~\textrm{then},
\end{displaymath}
\begin{displaymath}
|\Omega (G,\cdot )|=|\Theta (G,\cdot)||SA(G,\cdot )|.
\end{displaymath}
\end{myth}
{\bf Proof}\\
Let $A,B\in \Omega(G,\cdot )$. Then there exist $\alpha, \beta\in
SSYM (G,\cdot)$ and some $f_1,g_1,f_2,g_2\in H$ such that
\begin{displaymath}
A=(\alpha R_{g_1}^{-1},\alpha L_{f_1}^{-1},\alpha ),~B=(\beta
R_{g_2}^{-1},\beta L_{f_2}^{-1},\beta )\in AUT(G,\cdot ).
\end{displaymath}
Define a relation $\thicksim$ on $\Omega(G,\cdot )$ such that
\begin{displaymath}
A\thicksim B\Longleftrightarrow f_1=f_2~\textrm{and}~g_1=g_2.
\end{displaymath}
It is very easy to show that $\thicksim$ is an equivalence relation
on  $\Omega(G,\cdot )$. It can easily be seen that the equivalence
class $[A]$ of $A\in \Omega(G,\cdot )$ is the inverse image of the
mapping
\begin{displaymath}
\Psi~:~\Omega (G,\cdot )\longrightarrow \Theta (G,\cdot
)~\textrm{defined as}~\Psi~:~(\alpha R_{g_1}^{-1},\alpha
L_{f_1}^{-1},\alpha )\mapsto (f,g).
\end{displaymath}
If $A,B\in \Omega(G,\cdot )$ then $\Psi (A)=\Psi (B)$ if and only if
$(f_1,g_1)=(f_2,g_2)$ so, $f_1=f_2$ and $g_1=g_2$. Thus, since
$\Omega (G,\cdot )\le AUT(G,\cdot)$ by Theorem~\ref{1:15}, then
$AB^{-1}=(\alpha R_{g_1}^{-1},\alpha L_{f_1}^{-1},\alpha )(\beta
R_{g_2}^{-1},\beta L_{f_2}^{-1},\beta )^{-1}=(\alpha
R_{g_1}^{-1}R_{g_2}\beta^{-1},\alpha
L_{f_1}^{-1}L_{f_2}\beta^{-1},\alpha\beta^{-1} )=(\alpha
\beta^{-1},\alpha \beta^{-1},\alpha\beta^{-1} )\in AUT(G,\cdot
)\Leftrightarrow\alpha\beta^{-1}\in SA(G,\cdot )$. So,
\begin{displaymath}
A\thicksim B\Longleftrightarrow \alpha\beta^{-1}\in SA(G,\cdot
)~\textrm{and}~(f_1,g_1)=(f_2,g_2).
\end{displaymath}
$\therefore |[A]|=|SA(G,\cdot )|$. But each $A=(\alpha
R_g^{-1},\alpha L_f^{-1},\alpha )\in \Omega(G,\cdot )$ is determined
by some $f,g\in H$. So since the set $\Big\{[A]~:~A\in\Omega
(G,\cdot )\Big\}$ of all equivalence classes partitions
$\Omega(G,\cdot )$ by the fundamental theorem of equivalence
Relation,
\begin{displaymath}
|\Omega (G,\cdot )|=\sum_{f,g\in H}|[A]|=\sum_{f,g\in H}|SA(G,\cdot
)|=|\Theta (G,\cdot)||SA(G,\cdot )|.
\end{displaymath}
$\therefore |\Omega (G,\cdot )|=|\Theta (G,\cdot)||SA(G,\cdot )|$.

\begin{myth}\label{1:20}
Let $(G,\cdot )$ be a finite Smarandache loop with a finite
S-subgroup $H$. $(G,\cdot )$ is S-isomorphic to all its S-loop
S-isotopes if and only if
\begin{displaymath}
|(H,\cdot )|^2|SA(G,\cdot )|=|SBS(G,\cdot)||N_\mu (G,\cdot )|.
\end{displaymath}
\end{myth}
{\bf Proof}\\
As shown in [Corollary~5.2, \cite{sma8}], an S-loop is S-isomorphic
to all its S-loop S-isotopes if and only if it is S-isomorphic to
all its Smarandache $f,g$ principal isotopes. This will happen if
and only if $H\times H=\Theta (G,\cdot )$ where $\Theta (G,\cdot )$
is as defined in Theorem~\ref{1:19}.

Since $\Theta (G,\cdot )\subseteq H\times H$ then it is easy to see
that for a finite Smarandache loop with a finite S-subgroup $H$,
$H\times H=\Theta (G,\cdot )$ if and only if $|H|^2=|\Theta
(G,\cdot)|$. So the proof is complete by Theorem~\ref{1:18} and
Theorem~\ref{1:19}.

\begin{mycor}\label{1:21}
Let $(G,\cdot )$ be a finite Smarandache loop with a finite
S-subgroup $H$. $(G,\cdot )$ is a GS-loop if and only if
\begin{displaymath}
|(H,\cdot )|^2|SA(G,\cdot )|=|SBS(G,\cdot)||N_\mu (G,\cdot )|.
\end{displaymath}
\end{mycor}
{\bf Proof}\\
This follows by the definition of a GS-loop and Theorem~\ref{1:20}.

\begin{mylem}\label{1:22}
Let $(G,\cdot )$ be a finite GS-loop with a finite S-subgroup $H$
and a middle nucleus $N_\mu (G,\cdot )$ .
\begin{displaymath}
|(H,\cdot )|=|N_\mu (G,\cdot )|\Longleftrightarrow |(H,\cdot
)|=\frac{|SBS(G,\cdot)|}{|SA(G,\cdot )|}.
\end{displaymath}
\end{mylem}
{\bf Proof}\\
From Corollary~\ref{1:21},
\begin{displaymath}
|(H,\cdot )|^2|SA(G,\cdot )|=|SBS(G,\cdot)||N_\mu (G,\cdot )|.
\end{displaymath}
\begin{description}
\item[Necessity]
If $|(H,\cdot )|=|N_\mu (G,\cdot )|$, then
\begin{displaymath}
|(H,\cdot )||SA(G,\cdot )|=|SBS(G,\cdot)|\Longrightarrow |(H,\cdot
)|=\frac{|SBS(G,\cdot)|}{|SA(G,\cdot )|}.
\end{displaymath}
\item[Sufficiency]
If $|(H,\cdot )|=\frac{|SBS(G,\cdot)|}{|SA(G,\cdot )|}$ then,
$|(H,\cdot )||SA(G,\cdot )|=|SBS(G,\cdot)|$. Hence, multiplying both
sides by $|(H,\cdot )|$,
\begin{displaymath}
|(H,\cdot )|^2|SA(G,\cdot )|=|SBS(G,\cdot)||(H,\cdot )|.
\end{displaymath}
So that
\begin{displaymath}
|SBS(G,\cdot)||N_\mu (G,\cdot )|=|SBS(G,\cdot)||(H,\cdot
)|\Longrightarrow |(H,\cdot )|=|N_\mu (G,\cdot )|.
\end{displaymath}
\end{description}

\begin{mycor}\label{1:23}
Let $(G,\cdot )$ be a finite GS-loop with a finite S-subgroup $H$.
If $|N_\mu (G,\cdot )|\gneqq 1$, then,
\begin{displaymath}
|(H,\cdot )|=\frac{|SBS(G,\cdot)|}{|SA(G,\cdot
)|}.~\textrm{Hence},~|(G,\cdot )|=
\frac{n|SBS(G,\cdot)|}{|SA(G,\cdot )|}~\textrm{for some $n\gneqq
1$}.
\end{displaymath}
\end{mycor}
{\bf Proof}\\
By hypothesis, $\{e\}\ne H\ne G$. In a loop, $N_\mu (G,\cdot )$ is a
subgroup, hence if $|N_\mu (G,\cdot )|\gneqq 1$, then, we can take
$(H,\cdot )=N_\mu (G,\cdot )$ so that $|(H,\cdot )|=|N_\mu (G,\cdot
)|$. Thus by Lemma~\ref{1:22}, $|(H,\cdot
)|=\frac{|SBS(G,\cdot)|}{|SA(G,\cdot )|}$.

As shown in [Section~1.3, \cite{phd42}], a loop $L$ obeys the
Lagrange's theorem relative to a subloop $H$ if and only if
$H(hx)=Hx$ for all $x\in L$ and for all $h\in H$. This condition is
obeyed by $N_\mu (G,\cdot )$, hence
\begin{displaymath}
|(H,\cdot )|\Big||(G,\cdot)|\Longrightarrow
\frac{|SBS(G,\cdot)|}{|SA(G,\cdot
)|}\Bigg||(G,\cdot)|\Longrightarrow
\end{displaymath}
there exists $n\in \mathbb{N}$ such that
\begin{displaymath}
|(G,\cdot )|= \frac{n|SBS(G,\cdot)|}{|SA(G,\cdot )|}.
\end{displaymath}
But if $n=1$, then $|(G,\cdot )|=|(H,\cdot )|\Longrightarrow
(G,\cdot )=(H,\cdot )$ hence $(G,\cdot )$ is a group which is a
contradiction to the fact that $(G,\cdot )$ is an S-loop.
\begin{displaymath}
\therefore~|(G,\cdot )|=\frac{n|SBS(G,\cdot)|}{|SA(G,\cdot
)|}~\textrm{for some natural numbers $n\gneqq 1$}.
\end{displaymath}

\end{document}